\documentclass[10pt,a4paper]{article}
\usepackage{amsmath,amsfonts,amsthm,amsopn}

\newcommand{\eps}{\varepsilon}
\newcommand{\R}{\mathbb{R}}

\newcommand{\RN}{{\mathbb{R}^N}}

\newcommand{\de}{\partial}

\DeclareMathOperator{\cat}{cat}
\DeclareMathOperator{\dv}{div}
\DeclareMathOperator{\dist}{dist}

\renewcommand{\d }{\delta }

\newcommand{\n }{\nabla }

\renewcommand{\O}{\Omega}
\renewcommand{\OE}{\Omega_\varepsilon}

\newcommand{\G}{\Gamma}

\newcommand{\A}{{\cal A}}

\renewcommand{\P}{\mathcal{P}}

\newtheorem{theorem}{Theorem}[section]
\newtheorem{lemma}[theorem]{Lemma}

\newtheorem{proposition}[theorem]{Proposition}
\newtheorem{remark}[theorem]{Remark}

\begin{document}

\date{ }

\title{\textbf{Interior spikes of a singularly perturbed Neumann problem with potentials}}
\author{
Alessio Pomponio\thanks{Both authors supported by MIUR, national project \textit{Variational methods and nonlinear differential equations.}} \\ SISSA, via Beirut 2/4, I-34014 Trieste. \\ {\tt pomponio@sissa.it} \and Simone Secchi \\ Dipartimento di Matematica ``F. Enriques" \\
Universit\`a degli Studi di Milano \\ via C. Saldini 50, I-20133 Milano. \\
{\tt secchi@mat.unimi.it}}
\maketitle

{\small
\noindent \textbf{Abstract} \quad In this paper we prove that a singularly perturbed Neumann 
problem with potentials 
admits the existence of interior spikes concentrating in maxima and minima of an auxiliary functional 
depending only on the potentials.}

{\small
\noindent \textbf{Keywords} \quad Interior spikes, singularly perturbed Neumann problem.}

\section{Introduction}

In this paper we study the existence of interior spikes of the
following problem:
\begin{equation}\label{EQe}
\left\{
\begin{array}[c]{ll}
-\eps^2 \dv \left(J(x)\nabla u\right)+V(x)u=u^p & \text{in }\O,
\\
\frac{\de u}{\de \nu}=0 & \text{on }\de\O,
\end{array}
\right.
\end{equation}
where $\O$ is a smooth bounded domain of $\RN$ with external normal $\nu$,
$N\ge 3$, $1<p<(N+2)/(N-2)$, $J\colon\RN \to \R$ and $V\colon\RN
\to \R$ are $C^2$ functions.

In \cite{P}, the first author, extending the classical results by Ni and Takagi, in \cite{NT, NT2},
proved that there exist solutions of
\eqref{EQe} that concentrate at maximum and minimum points of a
suitable auxiliary function defined on the boundary $\de \O$ and
depending only on $J$ and $V$. Here we study the existence of
solutions which concentrate in the interior of $\O$ and we will
show that the concentration occurs at maximum and minimum points
of the same auxiliary function introduced in \cite{P}, but now
defined in $\O$. We assume that the reader has familiarity with \cite{P}.

When $J\equiv 1$ and $V\equiv 1$, interior spikes have been found by Wei (see \cite{W}) showing that concentration occurs at local maxima of the
distance function $\dist(\cdot, \de \O)$.

On $J$ and $V$ we will do the following assumptions:
\begin{description}
\item[(J)] $J\in C^2 (\O, \R)$, $J$ and $D^2 J$ are bounded; moreover, 
$J(x)\ge C>0$ for all $x\in\O$;
\item[(V)] $V\in C^2 (\O, \R)$, $V$ and $D^2 V$ are bounded; moreover,
$V(x)\ge C>0$ for all $x\in\O$.
\end{description}

Let us introduce an auxiliary function which will play a crucial
r\^ole in the study of \eqref{EQe}. Let $\G\colon \O \to \R$ be the
function defined by:
\begin{equation}\label{eq:Gamma}
\G(Q)=V(Q)^{\frac{p+1}{p-1}- \frac{N}{2}} J(Q)^{\frac{N}{2}}.
\end{equation}
Let us observe that by {\bf (J)} and {\bf (V)}, $\G$ is well
defined. Our main result is:

\begin{theorem}\label{th1}
Let  $Q_0 \in \O$. Suppose {\bf (J)} and {\bf (V)} hold. There exists
$\eps_0>0$ such that if $0<\eps<\eps_0$, then \eqref{EQe}
possesses a solution $u_\eps$ concentrating at $Q_\eps$ with
$Q_\eps \to Q_0$, as $\eps \to 0$, provided that one of the two
following conditions holds:
\begin{description}
\item[$(a)$] $Q_0$ is a non-degenerate critical point of $\G$;
\item[$(b)$] $Q_0$ is an isolated local strict minimum or maximum of $\G$.
\end{description}
\end{theorem}

%%%%%%%%%%%%%%%%%%%%%%%%%%%%%%%%%%%%%%%%%%%%%%%%%%%%%%%%%%%%%%%%%%%%%%%%%%%%%%%%%%%

\begin{center}{\bf Notation}\end{center}
\begin{itemize}
\item If $u \colon \RN \to \R$ and $P\in \RN$, we set $u_P \equiv u(\cdot - P)$.

\item If $U^Q$ is the function defined in \eqref{eq:UQ}, when there is no misunderstanding,
we will often write $U$ instead of $U^Q$. Moreover if $P=Q/\eps$, then $U_P \equiv U^Q(\cdot -P)$.

\item If $\eps>0$, we set $\OE \equiv \O/\eps \equiv\{x\in \RN : \eps x \in \O\}$.

\item With $o_\eps(1)$ we denote a quantity which tends to zero as $\eps \to 0$.
\end{itemize}

%%%%%%%%%%%%%%%%%%%%%%%%%%%%%%%%%%%%%%%%%%%%%%%%%%%%%%%%%%%%%%%%%%%%%%%%%%%%%%%%%%%%%%%%%%%%%%%

\section{Preliminary lemmas and some estimates}

First of all we perform the change of variables $x \mapsto \eps x$ and so problem \eqref{EQe} becomes
\begin{equation}\label{EQ}
\begin{cases}
-\dv \left(J(\eps x)\nabla u\right)+V(\eps x)u=u^p
& \text{in }\O_\eps,
\\
\frac{\de u}{\de \nu}=0
& \text{on }\de\O_\eps,
\end{cases}
\end{equation}
where $\O_\eps=\eps^{-1}\O$. Of course if $u$ is a solution of \eqref{EQ}, then $u(\cdot/\eps)$
is a solution of \eqref{EQe}.

Solutions of (\ref{EQ}) are critical points $u\in H^1(\O_\eps)$ of
\[
f_\eps (u)=
\frac{1}{2}\int_{\OE} J(\eps x)|\nabla u|^2 dx+
\frac{1}{2}\int_{\OE} V(\eps x)u^2dx
-\frac{1}{p+1}\int_{\OE} |u|^{p+1}.
\]
We look for solutions of (\ref{EQ}) near a $U^Q$, the unique solution of the \textit{limiting problem}
\begin{equation}\label{eq:Q}
\begin{cases}
-J(Q)\varDelta  u+ V(Q) u= u^p &   \text{in }\R^{N},\\
u>0 &   \text{in }\R^{N},\\
u(0)=\max_{\R^N} u, &
\end{cases}
\end{equation}
for an appropriate choice of $Q \in \O$. It is easy to see that
\begin{equation}\label{eq:UQ}
U^Q(x)=V(Q)^{\frac{1}{p-1}}\,\bar U\left(x \sqrt{V(Q)/J(Q)}  \right),
\end{equation}
where $\bar U$ is the unique solution of
\[
\begin{cases}
-\varDelta  \bar U+ \bar U= \bar U^p &   \text{in }\R^{N},
\\
\bar U>0 &   \text{in }\R^{N},
\\
\bar U(0)=\max_{\R^{N}}\bar U,
\end{cases}
\]
which is radially symmetric and decays exponentially at infinity together with its first derivatives.

For the sake of brevity, we will often write $U$ instead of $U^Q$.
If $P=\eps^{-1}Q \in \OE$, we set $U_P \equiv U^Q(\cdot - P)$ and
\[
Z^\eps \equiv \{ U_P : P\in \OE \}.
\]

\begin{lemma}\label{lem:nf}
For all $Q \in \O$ and for all $\eps$ sufficiently small, if $P=Q/\eps \in \OE$, then
\begin{equation}\label{eq:nf}
\|\nabla f_\eps(U_P)\|=O(\eps).
\end{equation}
\end{lemma}

\begin{proof}
Repeating the calculations of \cite{P}, we get:
\begin{eqnarray*}
(\nabla f_\eps(U_P) \mid v) &=&
 \int_{\frac{\O -Q}{\eps}} \left[-J(Q) \varDelta U + V(Q) U - U^p \right] v_{-P}
+J(Q) \int_{\de \OE} \frac{\de U_P}{\de \nu} v
\\
&&+\int_{\frac{\O -Q}{\eps}} (J(\eps x +Q)-J(Q)) \nabla U \cdot \n v_{-P}
+\int_{\frac{\O -Q}{\eps}} (V(\eps x +Q)-V(Q)) U v_{-P}.
\end{eqnarray*}
Hence, since $U\equiv U^Q$ is solution of \eqref{eq:Q}, we get
\begin{eqnarray}\label{eq:nf2}
(\nabla f_\eps(U_P) \mid v) &=&
J(Q) \int_{\de \OE} \frac{\de U_P}{\de \nu} v
+\int_{\frac{\O -Q}{\eps}} (J(\eps x +Q)-J(Q)) \nabla U \cdot \n v_{-P} \nonumber
\\
&&+\int_{\frac{\O -Q}{\eps}} (V(\eps x +Q)-V(Q)) U v_{-P}.
\end{eqnarray}
Let us estimate the first of these three terms:
\begin{gather*}
\Big| J(Q) \int_{\de \OE} \frac{\de U_P}{\de \nu} v \Big|
\le C \|v\|_{L^2(\de \OE)}
\Big( \int_{\de \OE}\Big| \frac{\de U_P}{\de \nu} \Big|^2 \Big)^{1/2}.
\end{gather*}
First of all, we observe that there exist $\eps_0>0$ and $C>0$ such that, for all
$\eps \in(0,\eps_0)$ and for all $v\in H^1(\OE)$, we have
\[
\|v\|_{L^2(\de \OE)}\le C \|v\|_{H^1(\OE)}.
\]
Using the exponential decay of $U$ and its derivatives, it is easy to see that
\begin{equation}\label{eq:deU}
\int_{\de \OE}\left|\frac{\de U_P}{\de \nu}\right|^2
=\int_{\de \left( \frac{\O-Q}{\eps}\right) }\left|\frac{\de U}{\de \nu}\right|^2
=o(\eps)
\end{equation}

Arguing as in \cite{P}, one can prove that:
\begin{equation}\label{eq:resto}
\int_{\frac{\O -Q}{\eps}} (J(\eps x +Q)-J(Q)) \nabla U  \cdot \n v_{-P}
+\int_{\frac{\O -Q}{\eps}} (V(\eps x +Q)-V(Q)) U v_{-P}
=O(\eps)\|v\|.
\end{equation}
Now the conclusion follows immediately from \eqref{eq:nf2}, \eqref{eq:deU} and \eqref{eq:resto}.
\end{proof}

We here present some useful estimates that will be used in the sequel.

\begin{proposition}\label{lemma1.2}
Let $P=Q/\eps \in \O_\eps$. Then we have:
\begin{equation*}%\label{eq:1.4}
\int_{\O_\eps} U_P^{p+1} =
\int_{\RN}\left(U^Q\right)^{p+1}
+o(\eps)
\end{equation*}

\begin{equation*}%\label{eq:1.5}
\int_{\de \O_\eps} \frac{\de  U_P}{\de \nu} U_P
=o(\eps),
\end{equation*}

\begin{equation*}%\label{eq:1.6}
J(Q) \int_{\O_\eps} |\n U_P|^2
+V(Q) \int_{\O_\eps} U_P^2
=\int_{\RN}\left(U^Q\right)^{p+1}
+o(\eps),
\end{equation*}

\begin{equation}\label{eq:J}
\int_{\O_\eps} J(\eps x)|\nabla U_P|^2=
J(Q) \int_{\OE}|\nabla U_P|^2
+\eps \int_{\RN}J'(Q)[x] |\nabla U^Q|^2
+o(\eps),
\end{equation}

\begin{equation}\label{eq:V}
\int_{\O_\eps}  V(\eps x)  U_P^2=
V(Q) \int_{\OE}  U_P^2
+\eps \int_{\RN} V'(Q)[x] \left(U^Q\right)^2
+o(\eps).
\end{equation}
\end{proposition}

\begin{proof}
Let us prove the first formula. We have:
\[
\int_{\OE} U_P^{p+1}=
\int_{\frac{\O-Q}{\eps}}\left(U^Q\right)^{p+1}=
\int_{\RN}\left(U^Q\right)^{p+1}
-\int_{\RN \setminus \frac{\O-Q}{\eps}}\left(U^Q\right)^{p+1}.
\]
Using the exponential of $U^Q$, it is easy to see that
\[
\int\limits_{\RN \setminus \frac{\O-Q}{\eps}}\left(U^Q\right)^{p+1}
\le \int\limits_{\RN \setminus B_{1/\sqrt\eps}}\left(U^Q\right)^{p+1}
= C \int_{\frac{1}{\sqrt\eps}}^\infty r^{N-1} (U^Q (r))^{p+1}\, dr
=o(\eps).
\]
Using the exponential of $U^Q$, also the second formula can be proved in a similar way. 
The other equations can be proved as in \cite{P}.
\end{proof}

%%%%%%%%%%%%%%%%%%%%%%%%%%%%%%%%%%%%%%%%%%%%%%%%%%%%%%%%%%%%%%%%%%%%%%%%%%%%%%%%%%%%%%%%%%%%%

\section{The finite dimensional reduction}

In this section we perform a finite dimensional reduction, following some ideas introduced in~\cite{AMS}.
The symbol $T_{U_P}Z^\eps$ denotes the tangent space to $Z^\eps$ at $U_P$.
 Let $L_{\eps,Q}:(T_{U_P}Z^\eps)^\perp\to
(T_{U_P}Z^\eps)^\perp$ denote the operator defined by setting
$(L_{\eps,Q}v \mid w)= D^2 f_\eps(U_P)[v,w]$.

\begin{lemma}\label{lem:inv}
There exists $C>0$ such that for $\eps$ small enough
one has that
\begin{equation}\label{eq:inv}
\|L_{\eps,Q}v \| \ge C \|v\|,\qquad \forall\; v\in(T_{U_P}Z^{\eps})^{\perp}.
\end{equation}
\end{lemma}

\begin{proof}
The proof of \eqref{eq:inv} is completely analogous to that of equation (21) in \cite{P}, so we omit the details.
\end{proof}

\begin{lemma}\label{lem:w}
For $\eps>0$ small enough, there exists a unique
$w=w(\eps, Q)\in (T_{U_P} Z^\eps)^{\perp}$ such that
$\nabla f_\eps (U_P + w)\in T_{U_P} Z$.
Such a $w(\eps,Q)$ is of class $C^{2}$, resp.  $C^{1,p-1}$, with respect to $Q$, provided
that $p\ge 2$, resp. $1<p<2$.
Moreover, the functional $\A_\eps (Q)=f_\eps (U_{Q/\eps} +w(\eps,Q))$ has
the same regularity of $w$ and satisfies:
\[
\nabla \A_\eps(Q_0)=0
\quad \Longleftrightarrow \quad
\nabla f_\eps\left(U_{Q_0/\eps}+w(\eps,Q_0)\right)=0.
\]
\end{lemma}

\begin{proof}
Let $\P=\P_{\eps, Q}$ denote the projection onto $(T_{U_P} Z^\eps)^\perp$. We want
to find a solution $w\in (T_{U_P} Z^\eps)^{\perp}$ of the equation
$\P\nabla f_\eps(U_P +w)=0$.  One has that $\nabla f_\eps(U_P+w)=
\nabla f_\eps (U_P)+D^2 f_\eps(U_P)[w]+R(U_P,w)$ with $\|R(U_P,w)\|=o(\|w\|)$, uniformly
with respect to $U_P$. Therefore, our equation is:
\begin{equation*}
L_{\eps,Q}w + \P\nabla f_\eps (U_P)+\P R(U_P,w)=0.
\end{equation*}
According to Lemma \ref{lem:inv}, this is equivalent to
\[
w = N_{\eps,Q}(w), \quad \mbox{where}\quad
N_{\eps,Q}(w)=-L_{\eps,Q}\left( \P \nabla f_\eps (U_P)+\P R(U_P,w)\right).
\]
By \eqref{eq:nf} it follows that
\begin{equation}\label{eq:N}
\|N_{\eps,Q}(w)\| = O(\eps) + o(\|w\|).
\end{equation}
Then one readily checks that $N_{\eps,Q}$ is a contraction on some ball in
$(T_{U_P} Z^\eps)^{\perp}$
provided that $\eps>0$ is small enough.
Then there exists a unique $w$ such that $w=N_{\eps,Q}(w)$.  Given $\eps>0$ small,
we can apply the Implicit
Function Theorem to the map $(Q,w)\mapsto \P \nabla f_\eps (U_P + w)$.
Then, in particular, the function $w(\eps,Q)$ turns out to be of class
$C^1$ with respect to $Q$.  Finally, it is a standard argument, see
\cite{AMS}, to check that the critical points of $\A_\eps(Q)=f_\eps (U_P+w)$
give rise to critical points of $f_\eps$.
\end{proof}

\begin{remark}\label{rem:w}
From (\ref{eq:N}) it follows immediately that:
\begin{equation}\label{eq:w}
\|w\|=O(\eps).
\end{equation}
Moreover repeating the arguments of \cite{P}, if $\gamma=\min\{1,p-1\}$, then, for $i=1, \ldots, N$,
we infer that
\begin{equation*}
\|\de_{P_i} w\|=O(\eps^\gamma).
\end{equation*}
\end{remark}

\section{The finite dimensional functional}

\begin{theorem}\label{th:sviluppo}
Let $Q \in \O$ and $P=Q/\eps \in \OE$. Suppose {\bf (J)} and {\bf (V)}. Then,
for $\eps$ sufficiently small, we get:
\begin{equation*}
\A_\eps (Q) =
f_\eps(U_P + w(\eps,Q))
= c_0 \G(Q)
+\frac{\eps}{2}\int_{\RN}J'(Q)[x] |\nabla U|^2
+\frac{\eps}{2}\int_{\RN}V'(Q)[x] U^2
+o(\eps),
\end{equation*}
where $\G$ is the auxiliary function introduced in \eqref{eq:Gamma} and
\[
c_0\equiv\left(\frac{1}{2}-\frac{1}{p+1}\right)\int_{\RN} \bar U^{p+1}.
\]
Moreover, for all $i=1,\ldots,N$, we get:
\begin{equation}\label{eq:DA}
\de_{P_i} \A_\eps (Q)= \eps c_0 \de_{Q_i} \G(Q)+o(\eps),
\end{equation}
and hence
\begin{equation}\label{eq:C1}
\|\A_\eps - c_0 \G \|_{C^1(\O)}=O(\eps).
\end{equation}

\end{theorem}

\begin{proof}
In the sequel, to be short, we will often write $w$ instead of $w(\eps,Q)$.
It is always understood that $\eps$ is taken in such a
way that all the results discussed previously hold.

First of all, reasoning as in the proofs of \eqref{eq:J} and \eqref{eq:V}
and by \eqref{eq:w}, we can observe that
\begin{eqnarray}
\int_{\OE} J(\eps x)\nabla U_P \cdot \nabla w
&=&
J(Q) \int_{\OE} \nabla U_P \cdot \nabla w + o(\eps),     \label{eq:J2}
\\
\int_{\OE} V(\eps x)U_P \,w
&=&
V(Q) \int_{\OE} U_P \,w +o(\eps).                        \label{eq:V2}
\end{eqnarray}
Recalling \eqref{eq:w}, we have:
\begin{gather*}
\A_\eps (Q)
=\frac{1}{2}\int_{\OE} J(\eps x) |\nabla (U_P +w)|^2
+\frac{1}{2}\int_{\OE} V(\eps x)(U_P+w)^2
-\frac{1}{p+1}\int_{\OE} (U_P+w)^{p+1}
\\
=\frac{1}{2}\int_{\OE} J(\eps x) |\nabla U_P|^2
+\frac{1}{2}\int_{\OE} V(\eps x)U_P^2
-\frac{1}{2}\int_{\OE} U_P^{p+1}
+\int_{\OE} J(\eps x)\nabla U_P \cdot \nabla w
+\int_{\OE} V(\eps x)U_P \,w
\\
 -\int_{\OE} U_P^p \,w
+\left(\frac{1}{2}-\frac{1}{p+1}\right) \int_{\OE} U_P^{p+1}
-\frac{1}{p+1}\int_{\OE}\left[ (U_P+w)^{p+1} -U_P^{p+1}-(p+1) U_P^p \,w\right]
+o(\eps)=
\end{gather*}
[by Proposition \ref{lemma1.2}, \eqref{eq:J2} and \eqref{eq:V2}
and with our notations]
\begin{eqnarray*}
&=& \frac{1}{2} \int_{\RN} U^{p+1}
+\frac{\eps}{2}\int_{\RN}J'(Q)[x] |\nabla U|^2
+\frac{\eps}{2}\int_{\RN}V'(Q)[x] U^2
-\frac{1}{2}\int_{\RN} U^{p+1}
\\
&&+J(Q) \int_{\OE}  \nabla U_P \cdot \nabla w
+V(Q) \int_{\OE}  U_P \,w
-\int_{\OE} U_P^p \,w
+\left(\frac{1}{2}-\frac{1}{p+1}\right) \int_{\RN} U^{p+1}
+o(\eps).
\end{eqnarray*}
Since that $U$ is solution of \eqref{eq:Q}, we infer
\begin{gather*}
J(Q) \int_{\OE}  \nabla U_P \cdot \nabla w
+V(Q) \int_{\OE}  U_P \,w
-\int_{\OE} U_P^p \,w
\\
= \int_{\OE} \left[ -J(Q) \varDelta U_P + V(Q)U_P - U_P^p \right] w
+J(Q)\int_{\de \OE} \frac{\de U_P}{\de \nu} w
= J(Q)\int_{\de \OE} \frac{\de U_P}{\de \nu} w
=o(\eps).
\end{gather*}
By these considerations we can say that
\[
\A_\eps (Q) = \Big( \frac{1}{2}-\frac{1}{p+1} \Big) \int_{\RN} U^{p+1}
+\frac{\eps}{2}\int_{\RN}J'(Q)[x] |\nabla U|^2
+\frac{\eps}{2}\int_{\RN}V'(Q)[x] U^2
+o(\eps).
\]
Now the conclusion of the first part of the theorem follows observing that, since by \eqref{eq:UQ}
\[
U^Q (x)=V(Q)^{\frac{1}{p-1}}\,\bar U\left(x \sqrt{\tfrac{V(Q)}{J(Q)}}  \right),
\]
then
\[
\int_{\RN} U^{p+1}= V(Q)^{\frac{p+1}{p-1}-\frac{N}{2}}J(Q)^{\frac{N}{2}} \int_{\RN} \bar U^{p+1}.
\]
The estimate \eqref{eq:DA} on the derivatives of $\A_\eps$ follows easily by repeating the arguments 
of \cite{P}.
\end{proof}

%%%%%%%%%%%%%%%%%%%%%%%%%%%%%%%%%%%%%%%%%%%%%%%%%%%%%%%%%%%%%%%%%%%%%%%%%%%%%%%

\section{Proof of Theorem \ref{th1}}

In this section we will state and prove two multiplicity results for \eqref{EQe}.
Theorem \ref{th1} will follow as a particular case.

\begin{theorem}\label{th:lc}
Let {\bf (J)} and {\bf (V)} hold and suppose $\G$ has a nondegenerate smooth manifold of critical points
$M\subset \O$.
There exists $\eps_0>0$ such that if $0<\eps<\eps_0$, then \eqref{EQe} has at least $l(M)$ solutions that
concentrate near points of $M$. Here $l(M)$ denotes the {\sl cup long} of $M$ (for a precise
definition, see \cite{C}).
\end{theorem}

\begin{proof}
Fix a $\delta$-neighborhood $M_\delta$ of $M$ such that
the only critical points of $\G$ in $M_\delta$ are those in $M$. We will take $U=M_\delta$. 
For $\eps$ sufficiently small, by \eqref{eq:C1} and Theorem 6.4 in Chapter II of \cite{C},
$\A_\eps$ possesses at least $l(M)$ critical points, which are solutions of \eqref{EQ}
by Lemma \ref{lem:w}. Let $Q_\eps \in M$ be one of these critical points, then
$u_\eps^{Q_\eps}=U_{Q_\eps/\eps}+w(\eps, Q_\eps)$ is a solution of \eqref{EQ}. Therefore
\[
u_\eps^{Q_\eps}(x/\eps)\simeq U_{Q_\eps/\eps}(x/\eps)= U^{Q_\eps}\left(\frac{x-Q_\eps}{\eps} \right)
\]
is a solution of \eqref{EQe}.
\end{proof}

Moreover, when we deal with local minima (resp. maxima) of $\G$, the
preceding results can be improved because the number of positive solutions of \eqref{EQe}
can be estimated by means of the category and $M$ does not need to be a manifold. 
We will give only the statement of the theorem; for the proof, see \cite{P}.

\begin{theorem}\label{th:cat}
Let {\bf (J)} and {\bf (V)} hold and suppose $\G$ has
a compact set $X \subset \O$ where $\G$ achieves a strict local minimum (resp. maximum),
in the sense that there exist $\delta>0$ and a $\d$-neighborhood $X_\d \subset \O$ of $X$ such that
\[
b\equiv \inf\{\G(Q) : Q\in \partial X_{\d}\}> a\equiv \G_{|_X}, \quad
\left({\rm resp. }\; \sup\{\G(Q) : Q\in \partial X_{\d}\}< \G_{|_X} \right).
\]
Then there exists $\eps_0>0$ such that \eqref{EQe}
has at least $\cat(X,X_\d)$ solutions that concentrate near points of $X_{\d}$, provided
$\eps \in (0,\eps_0)$. Here $\cat(X,X_\d)$ denotes the Lusternik-Schnirelman category of $X$ with
respect to $X_\d$.
\end{theorem}

\begin{remark}
Let us observe that part (a) of Theorem \ref{th1} is a particular case of Theorem \ref{th:lc} while
part (b) of Theorem \ref{th1} is a particular case of Theorem \ref{th:cat}.
\end{remark}


\begin{thebibliography}{99}



\bibitem{AMS}
A. Ambrosetti, A. Malchiodi, S. Secchi, {\it Multiplicity results
for some nonlinear Schr\"odinger equations with potentials}, Arch.
Rational Mech. Anal., {\bf 159}, (2001), 253--271.

\bibitem{C}
K. C. Chang, Infinite dimensional Morse theory and multiple
solutions problems, Birkh\"auser, 1993.

\bibitem{NT}
W. M. Ni, I. Takagi, {\it On the shape of the least-energy
solutions to a semilinear Neumann problem}, Comm. Pure Appl.
Math., {\bf 44}, (1991), 819--851.

\bibitem{NT2}
W. M. Ni, I. Takagi, {\it Locating the peaks of least-energy
solutions to a semilinear Neumann problem}, Duke Math. J., {\bf
70}, (1993), 247--281.

\bibitem{P}
A. Pomponio, {\it Singularly perturbed Neumann problems with
potentials}, Preprint SISSA, 87/2003/M, (2003). 
{\tt http://www.mat.uniroma3.it/AnalisiNonLineare/preprints/preprints.html}.

\bibitem{W}
J. Wei, {\it On the interior spike layer solutions to a singularly
perturbed Neumann problem}, Tohoku Math. J., {\bf 50}, (1998), 159--178.
\end{thebibliography}
\end{document}